\input amstex
\documentstyle{amsppt}
\NoRunningHeads
\TagsOnRight

\def\H{\Cal H}
\def\l{\langle}
\def\r{\rangle}
\def\wo{\widehat\otimes}
\topmatter
\title
Symmetric Extensions of Dirichlet Operators
\endtitle
\author
A. G. Us \\ Institute of Mathematics \\
Ukrainian National Academy of Science
\endauthor
\endtopmatter
\document

Classical Dirichlet forms and operators are studied intensively in
infinite-dimensi\-onal analysis and its applications (see [1] and references
therein). Recently there have been studied supersymmetric Dirichlet forms
and operators [2--3] because of numerous applications to quantum physics.

It this notice there is constructed and considered the extension of a
classical Dirichlet operator in the space of symmetric differential forms.

The work consists of 3 parts. In the first part there're recalled some
basic concepts and definitions concerning smooth measures on Hilbert spaces
and corresponding Dirichlet forms and operators. In the part 2 there's
constructed the extension of a classical Dirichlet operators and explicit
form is given.

In the third part the extension is considered as a differential operator in
the space of square integrated functions. There are given sufficient
conditions for its essential self-adjointness in onedimensional case.
Besides for the "supersymmetric part" of the operator conditions for
self-adjointness in general situation are found.

1. Let $\H$ be a real separable Hilbert space with a scalar product
$\l\cdot,\cdot\r$ and a norm $|\cdot|$ and $\Phi'\supset\H\supset\Phi$ be a
rigging of $\H$ by nuclear Frech\'et space $\Phi$ (densely and continuously
embedded) and its dual space $\Phi'$. The duality between $\Phi$ and
$\Phi'$, which is given by the scalar product in $\H$, will be also denoted
by $\l\cdot,\cdot\r$.

Denote by $\Cal FC_b^\infty(\Phi)$ the set of all smooth cylinder functions
on $\Phi'$ with all derivatives bounded. For $u\in\Cal FC_b^\infty(\Phi)$
and $\varphi\in\Phi$ $\nabla_\varphi u(\cdot)=\l\nabla u(\cdot),\varphi\r$
is a directional derivative of a $u(\cdot)$.

Let $\mu$ be a probability measure on $(\Phi',\Cal B(\Phi'))$ which is
quasiinvariant under translations by all elements of $\Phi$. We also
suppose that $\mu$ has a logarithmic derivative $\beta_\mu(\cdot)$
($\beta_\mu(\cdot):\Phi'\to\Phi'$ measurable mapping (see [4])) and
$\beta_\mu(\cdot)$ is weakly differentiable for $\mu$-a.e. $x\in\Phi'$.
Denote by $R_\mu(x)=-\beta_\mu'(x)$, $x\in\Phi'$ a family of self-adjoint
(possibly unbounded) operators in $\H$, such that $\Phi\subset
D(R_\mu(\cdot))$ $(\roman{mod}\mu)$. We shall say that the measure $\mu$ is
uniformly log-concave (ULC) if $\exists C>0$: $\forall\varphi\in\Phi$
$\l R_\mu(x)\varphi,\varphi\r\geqslant C|\varphi|^2$.

In our further considerations we assume that the ULC measure $\mu$
satisfies the following integrability conditions:
$$
\forall\varphi\in\Phi
\ \ \ \ \ \ \ \ \ \ \int_{\Phi'}\beta_{\mu,\varphi}^2(x)\roman d\mu(x)
<\infty, \tag1
$$
$$
\forall\varphi\in\Phi
\ \ \ \ \ \ \ \ \ \ \int_{\Phi'}|R_\mu(x)\varphi|^2\roman d\mu(x)
<\infty, \tag2
$$
where $\beta_{\mu,\varphi}(x):=\l\beta_\mu(x),\varphi\r$, $x\in\Phi'$.

On the domain $D(H_\mu)=\Cal FC_b^\infty(\Phi)\ $ define the differential
operator $H_\mu$ in
$\roman L_2(\mu)\equiv\roman L_2(\Phi',\Cal B(\Phi'),\mu)$:
$$
(H_\mu u)(x)=-\Delta u(x)-\l\beta_\mu(x),\nabla u(x)\r,\ x\in\Phi',
$$
where $\Delta u(\cdot)=Tr_\H u''(\cdot)$. Taking into account \thetag1 we
can see that the operator $H_\mu$ is well defined in $\roman L_2(\mu)$.
This operator is called the Dirichlet operator corresponding to the measure
$\mu$.

For $u,v\in\Cal FC_b^\infty(\Phi)$ define a positive symmetric form
$$
\Cal E_\mu(u,v):=(H_\mu u,v)_{\roman L_2(\mu)}
=\int_{\Phi'}\l\nabla u(x),\nabla\overline v(x)\r\roman d\mu(x).
$$

The form $\Cal E_\mu$ is obviously closable and its closure is a classical
Dirichlet form.

2. Consider the space
$$
\Gamma_\mu(\H):=\roman L_2(\mu)\otimes\Gamma(\H)
=\underset{n=0}\to{\overset\infty\to\oplus}\roman L_2(\mu)\otimes
\Gamma^n(\H)=\underset{n=0}\to{\overset\infty\to\oplus}\Gamma_\mu^n(\H),
$$
where $\Gamma(\H)$ is a symmetric (bosonic) Fock space and $\Gamma^n(\H)$,
$n\geqslant 0$ are $n$-particle subspaces ($\Gamma^0(\H)=\Bbb C$). The 
following
representation of $\Gamma_\mu(\H)$ is helpful:
$\Gamma_\mu(\H)=\roman L_2(\Phi'\to\Gamma(\H),\mu)$, for
$F\in\Gamma_\mu(\H)$ $\|F\|^2
=\int_{\Phi'}\|F(x)\|^2_{\Gamma(\H)}\roman d\mu(x)$.

For $n\geqslant 0$ the set $D_n(\Phi):=l.s.\{u_n(x)=u(x)\varphi_1\wo\dots
\wo\varphi_n:\ u\in\Cal FC_b^\infty(\Phi),\ \varphi_k\in\Phi\}$ is dense in
$\Gamma_\mu^n(\H)$. Then the set $D(\Phi)=\{(u_0,u_1,\dots,u_N,0,\dots,0):
\ u_k\in D_k(\Phi)\}$ is dense in $\Gamma_\mu(\H)$.

For any $n\geqslant 0$ introduce the linear operator
$\delta_n:\Gamma_\mu^n(\H)\to\Gamma_\mu^{n+1}(\H)$ with the domain
$D_n(\Phi)$:
$$
(\delta_nu_n)(x)=\sqrt{n+1}\nabla u(x)\wo\varphi_1\wo\dots\wo\varphi_n
$$
Using the operator $\delta_n$, $n\geqslant 0$ on the domain $D(\Phi)$ define a
linear operator in $\Gamma_\mu(\H)$:
$$
(\delta u)_n=\delta_{n-1}u_{n-1},\ (\delta u)_0=0 \tag3
$$
A direct computation yields:
$$
(\delta_n^*u_{n+1})(x)
=-\frac1{\sqrt{n+1}}\sum_{j=1}^{n+1}(\nabla_{\varphi_j}u(x)
+\beta_{\mu,\varphi_j}(x)u(x))\varphi_1\wo\dots\widehat j\dots\wo\varphi_{n+1}.
$$
Analogously \thetag3 define the operator $\delta^*$ in $\Gamma_\mu(\H)$
with the domain $D(\Phi)$:
$$
(\delta^*u)_n=\delta_n^*u_{n+1}.
$$

Finally, on the domain $D(\Phi)$ define a linear operator in
$\Gamma_\mu(\H)$
$$
\Delta_\mu=\delta^*\delta+\delta\delta^*.
$$
Obviously $\Delta_\mu$ is a positive symmetric operator.

\proclaim{Theorem 1}
On $D(\Phi)$
$$
\Delta_\mu=H_\mu\otimes\bold 1+\bold 1\otimes\roman d\Gamma(R_\mu(\cdot))
+\Bbb A_\mu, \tag4
$$
where $\roman d\Gamma(R_\mu(\cdot))$ is the second quantization of
$R_\mu(\cdot)$ and $\Bbb A_\mu$ is a positive symmetric operator. Its
explicit form will be given below. As the conditions \thetag1 and \thetag2
are fulfilled all the operators in \thetag4 are well defined in
$\Gamma_\mu(\H)$.
\endproclaim

The operator $\Delta_\mu$ is symmetric extension of Dirichlet operator
$H_\mu$.

Point out that in [3] the extension was the same as in \thetag4 just
without the operator $\Bbb A_\mu$.

3. Let us make a transition from $\Gamma_\mu(\H)$ to $\roman L_2(\mu)\otimes
\roman L_2(\gamma)$ using the Segal isomorphism, $\gamma$ is the standard
Gaussian measure on $\Phi'$. Wee shall obtain (see [5], ch.~6])
$$
\gather
S\Gamma_\mu(\H)=\roman L_2(\mu)\otimes\roman L_2(\gamma),\\
SD(\Phi)=D=l.s.\Cal FC_b^\infty(\Phi)\otimes\Cal P(\Phi'),\\
\boldsymbol\Delta_\mu=S\Delta_\mu S^{-1}=H_\mu\otimes\bold 1
+\bold 1\otimes H_{\gamma,R_\mu(\cdot)}+\bold A_\mu=H_{\mu,\gamma}
+\bold A_\mu,
\endgather
$$
where $S=\bold 1\otimes I$, $\Cal P(\Phi')$ is the set of all continuous
polynomials on $\Phi'$, $H_{\gamma,R_\mu(\cdot)}$ is the Dirichlet operator
of the measure $\gamma$ with coefficient operator $R_\mu(\cdot)$, see [5,
ch.~6], and the operator $\bold A_\mu$ is defined by its action on
functions $u(\cdot)p(\cdot)\in D$:
$$
\multline (\bold A_\mu up)(x,y)=Tr_Hu''(x)p''(y)+
\l p''(y)\beta_\mu(x),\nabla u(x)\r\\
-\l u''(x)y,\nabla p(y)\r-\l\beta_\mu(x),y\r\l\nabla u(x),\nabla p(y)\r,
\ \ x,y\in\Phi',
\endmultline
$$
as $\beta_\gamma(y)=-y$, $y\in\Phi'$.

We can see that $\bold A_\mu$ is the fourth order differential operator of
a complex structure. That's why one managed to find the condition for
essential self-adjointness of the operator $\boldsymbol\Delta_\mu$ just in
case, if $\Phi'=\H=\Phi=\Bbb R$. We have $\boldsymbol\Delta_\mu=H_\mu\otimes
\bold 1+\bold 1\otimes H_{\gamma,R_\mu(\cdot)} +2H_\mu\otimes H_\gamma$ on
$C_0^\infty(\Bbb R)\otimes\Cal P(\Bbb R)$. The following result is true
(see [6]).

\proclaim{Theorem 2}
Let $f(\cdot)$ be the density of the measure $\mu$. If $f\in C^1(\Bbb R)$
and $\forall x\in\Bbb R$ $f(x)>0$, the operator $\boldsymbol\Delta_\mu$ is
essentially self-adjoint on $C_0^\infty(\Bbb R)\otimes\Cal P(\Bbb R)$.
\endproclaim

In infinite-dimensional case because of the complex structure of the
operator $\bold A_\mu$ there're no satisfactory conditions for essential
self-adjointness of $\boldsymbol\Delta_\mu$. But it's possible to find such
conditions for the operator $H_{\mu,\gamma}$.

It appears that $H_{\mu,\gamma}$ may be considered as the Dirichlet
operator of a certain perturbed measure on
$(\Phi'\times\Phi',\Cal B(\Phi'\times\Phi'))$. That's why in order to prove
essential self-adjointness of $H_{\mu,\gamma}$ one can use one of numerous
theorem, concerning classical Dirichlet operators (see e.g. [7]).

Let $\Phi'\supset\H_-\supset\H$, $\mu(\H_-)=\gamma(\H_-)=1$. Denote by
$g(x,y):=\|R_\mu^{-1/2}(x)y\|_-$ and $h(x,y)
:=\|\big(R_\mu'(x)R_\mu^{-1/2}(x)y\big)R_\mu^{-1/2}(x)y\|_-$, $x,y\in\H_-$,
where $R_\mu'(\cdot)$ is the weak derivative of $R_\mu(\cdot)$. Now we are
able to formulate

\proclaim{Theorem 3}
Assume that $\beta_\mu(\cdot)\in C_{b,loc}^3(\H_-,\H_-)$,
$\|\beta_\mu(\cdot)\|_-\in\roman L_2(\mu)$, $\{g,h\}\subset
\roman L_2(\mu\times\gamma)$. Then the operator $H_{\mu,\gamma}$ is
essentially self-adjoint on $D$.
\endproclaim

\bigskip
{\smc A. G. Us, postgraduate student, Institute of Mathematics of
\break Ukrainian Academy of Sciences, Terestenskovskaya st., 3, Kiev, \break 
Ukraine
}

\enddocument